\documentclass[11pt]{amsart}
\usepackage{amsfonts}
\usepackage{amsmath}
\usepackage{amsthm}
\usepackage{amssymb}

\newcommand{\bea}{\begin{eqnarray}}
\newcommand{\eea}{\end{eqnarray}}
\newcommand{\cal}{\mathcal}
\newtheorem{thm}{Theorem}[section]
\newtheorem{prop}[thm]{Proposition}
\newtheorem{lem}[thm]{Lemma}
\newtheorem{cor}[thm]{Corollary}
\newtheorem{rem}[thm]{Remark}
\newtheorem{defn}[thm]{Definition}

\begin{document}

\begin{abstract}A characterization of the space of
symmetric Laurent polynomials of type $(BC)_n$ which vanish
on a certain set of submanifolds is given by using the
Koornwinder-Macdonald polynomials. A similar characterization was
given previously for symmetric polynomials of type $A_n$
by using the Macdonald polynomials.
We use a new method which exploits the duality relation. The method
simplifies a part of the proof in the $A_n$ case.
\end{abstract}
\title[Zeros of Symmetric Laurent Polynomials of Type $(BC)_n$]
{Zeros of Symmetric Laurent Polynomials of Type $(BC)_n$
and Koornwinder-Macdonald Polynomials Specialized at
$t^{k+1}q^{r-1}=1$}
\author{Masahiro Kasatani}
\address{Department of Mathematics, Graduate School of Science, 
Kyoto University, Kyoto 606-8502, Japan}
\email{kasatani@math.kyoto-u.ac.jp}

\maketitle

%%%%%%%%%%          Section intro          %%%%%%%%%%

\section{Introduction}
~
Let $k,r,n$ be positive integers. We assume that $n\geq k+1$ and $r\geq 2$.
In \cite{FJMM}, $n$-variable symmetric polynomials satisfying certain zero
conditions are characterized by using the Macdonald polynomials \cite{Mac}
specialized at
\begin{equation}
t^{k+1}q^{r-1}=1.\label{resonance old}
\end{equation}
To be precise, the paper \cite{FJMM} works in the following setting.
Denote by $m$ the greatest common divisor of $k+1$ and $r-1$.
Let $\omega$ be an $m$-th primitive root of unity. Then, the variety
given by $t^{\frac{k+1}m}q^{\frac{r-1}m}=\omega$ is an irreducible component
of (\ref{resonance old}). It is uniformized as follows. Let
$\omega_1\in\mathbb{C}$ be such that $\omega_1^{(r-1)/m}=\omega$. We consider
the specialization of $t,q$ in terms of the uniformization parameter $u$,
\begin{equation}\label{SPEC}
t=u^{(r-1)/m},q=\omega_1u^{-(k+1)/m}.
\end{equation}
The following result was obtained in \cite{FJMM}.
\begin{thm}\label{thm:old}
For a partition $\lambda=(\lambda_1,\ldots,\lambda_n)$ satisfying
\begin{equation}
\lambda_i-\lambda_{i+k}\geq r \quad (1\leq i\leq n-k),\label{adm}
\end{equation}
the Macdonald polynomial $P_\lambda\in\mathbb{C}(t,q)[x_1,\ldots,x_n]^{{\mathfrak S}_n}$
has no pole at $(\ref{SPEC})$,
and when it is specialized at $(\ref{SPEC})$, it vanishes
on the submanifold given by
\begin{equation}
x_{i}/x_{i+1}=tq^{s_i} \qquad\hbox{ for $1\leq i \leq k$}\label{WHEEL}
\end{equation}
for each choice of non-negative integers $s_i$ such that
$\sum_{i=1}^ks_i\leq r-1$.
Conversely, the space of symmetric polynomials
$P\in\mathbb{C}(u)[x_1,\ldots,x_n]^{{\mathfrak S}_n}$ satisfying the above condition
is spanned by the Macdonald polynomials $P_\lambda$ specialized
at $(\ref{SPEC})$ where $\lambda$ satisfies $(\ref{adm})$.
\end{thm}
The condition that a polynomial vanishes on the submanifold (\ref{WHEEL})
is called the wheel condition corresponding to the submanifold (\ref{WHEEL})
and a partition $\lambda$ satisfying the condition
(\ref{adm}) is called a $(k,r,n)$-admissible partition.
Note that if we set $s_{k+1}=r-1-\sum_{i=1}^ks_i$, it follows that
$x_{k+1}/x_{1}=tq^{s_{k+1}}$ from (\ref{WHEEL}) and (\ref{resonance old}) .

In this paper, we obtain a similar result in the case of $n$-variable symmetric
Laurent polynomials of type $(BC)_n$. Here we say a Laurent polynomial
in the variables $x_1,\ldots,x_n$ is of type  $(BC)_n$ if and only if it is
symmetric and invariant for the change of the variable $x_1$ to $x_1^{-1}$.
The original case in \cite{FJMM} corresponds to $A_n$.
We use the Koornwinder-Macdonald polynomials $P_\lambda$ of type
$(BC)_n$ \cite{K} in order to characterize the space of symmetric Laurent
polynomials of type  $(BC)_n$ satisfying the wheel conditions.
The Koornwinder-Macdonald polynomials depend on six parameters
$t,q,a,b,c,d$.

We set $W_n:={\mathfrak S}_n \ltimes (\mathbb{Z}_2)^n$.
Our main result is
\begin{thm}\label{NEW}
Let $\lambda=(\lambda_1,\ldots,\lambda_n)$ be a $(k,r,n)$-admissible partition.
Then, the Koornwinder-Macdonald polynomial
 $P_\lambda\in\mathbb{C}(t,q,a,b,c,d)[x_1^{\pm1},\ldots,x_n^{\pm1}]^{W_n}$
 has no pole at $(\ref{SPEC})$,  and when it is specialized at $(\ref{SPEC})$,
it satisfies the wheel conditions corresponding to $(\ref{WHEEL})$.
Conversely, the space of symmetric Laurent polynomials of type $(BC)_n$ in 
$\mathbb{C}(u,a,b,c,d)[x_1^{\pm1},\ldots,x_n^{\pm1}]^{W_n}$ satisfying the wheel
conditions is spanned by the Koornwinder-Macdonald polynomials
$P_\lambda$ specialized at $(\ref{SPEC})$ where $\lambda$ are
$(k,r,n)$-admissible partitions.
\end{thm}

Although the statement of Theorem \ref{NEW} is quite analogous to
that of Theorem \ref{thm:old}, our proof of Theorem \ref{NEW} is
different from that of Theorem \ref{thm:old} given in \cite{FJMM}.
In fact, our method gives an alternative proof simpler than
the one given in \cite{FJMM} for the $A_n$ case.
We use the duality relation for the Koornwinder-Macdonald
 polynomials $P_\lambda$.
In \cite{KMSV}, we obtain a further result
by an application of the method used in this paper.

Let us explain the duality relation and the method of our proof.
We denote by $P^*_\lambda$ the {\it dual} Koornwinder-Macdonald polynomial
\cite{vD} defined by dual parameters $t,q,a^*,b^*,c^*,d^*$:
\begin{eqnarray}
\begin{array}{l}
a^*=-a^{1/2}b^{1/2}c^{1/2}d^{1/2}q^{-1/2},\quad b^*=-a^{1/2}b^{1/2}c^{-1/2}d^{-1/2}q^{1/2},\\
c^*=-a^{1/2}b^{-1/2}c^{1/2}d^{-1/2}q^{1/2},\quad d^*=-a^{1/2}b^{-1/2}c^{-1/2}d^{1/2}q^{1/2}.
\end{array}\label{dual_param}
\end{eqnarray}
For any partition $\mu=(\mu_1,\ldots,\mu_n)$ and
 $f\in\mathbb{C}(t,q,a,b,c,d)[x_1^{\pm1},\ldots,x_n^{\pm1}]$,
we define specializations $u_\mu(f)$, $u^*_\mu(f)$ of $f$ by
\begin{eqnarray}
\begin{array}{l}
u_\mu(f):=f(t^{n-1}q^{\mu_1}a^*,t^{n-2}q^{\mu_2}a^*,\cdots,q^{\mu_n}a^*),\\
u^*_\mu(f):=f(t^{n-1}q^{\mu_1}a,t^{n-2}q^{\mu_2}a,\cdots,q^{\mu_n}a).
\end{array}\label{special_u}
\end{eqnarray}
In particular, we have
\begin{eqnarray*}
u_0(f)&=&f(t^{n-1}a^*,t^{n-2}a^*,\ldots,a^*),\\
u^*_0(f)&=&f(t^{n-1}a,t^{n-2}a,\ldots,a).
\end{eqnarray*}
The duality relations reads as
\begin{equation}\label{DUAL}
\frac{u^*_\mu(P_\lambda)}{u^*_0(P_\lambda)}=
\frac{u_\lambda(P^*_\mu)}{u_0(P^*_\mu)}.
\end{equation}

To prove the two statements, (i) $P_\lambda$ has no pole at $(\ref{SPEC})$,
and (ii) $P_\lambda$ specialized at $(\ref{SPEC})$ satisfies
the wheel conditions corresponding to $(\ref{WHEEL})$,
we use the duality relation with special choices of $\mu$.
Here, we explain only the latter assuming that the former is already proved.
The details of the proofs are given in the main body of the paper.

In order to study the values of $P_\lambda$ on the submanifold (\ref{WHEEL}),
we use (\ref{DUAL}) by choosing $\mu=(\mu_1,\ldots,\mu_n)$ in such a way that
\begin{eqnarray}
&\mu_i-\mu_{i+1}=s_i \quad\mbox{ for $i=1,\ldots,k$},\label{WHL}\\
&\mu_i-\mu_{i+1}> 2[\frac{n}{k+1}](r-1) \quad\mbox{ for $i=k+1,\ldots,n-1$}.\label{WLD}
\end{eqnarray}
{}From the definition of the dual polynomial $P^*_\mu$, it has
no pole at the specialization (\ref{SPEC}) if (\ref{WLD}) is valid.
Without specialization (\ref{SPEC}) we have an explicit formula for
$u^*_0(P_\lambda)$ and $u_0(P^*_\mu)$, and we can easily count the order
 of zeros (or poles) for them. Using (\ref{DUAL}), we can prove that
$u^*_\mu(P_\lambda)$ vanishes at (\ref{SPEC}). Since there exist enough
$\mu$'s satisfying the conditions  (\ref{WHL}) and (\ref{WLD}),
the Laurent polynomial $P_\lambda$ itself should vanish at (\ref{WHEEL}).

This much is the proof of the first half of Theorem \ref{NEW}. 
Let $J^{(k,r)}$ be the space of symmetric Laurent polynomials $P$ of type
$(BC)_n$ satisfying the wheel conditions, and for a positive integer $M$,
let $J^{(k,r)}_M$ be its subspace consisting of $P$ such that
the degree of $P$ in each variable $x_i$ is less than $M$.
Because of the invariance for $x_i\leftrightarrow x_i^{-1}$,
the dimension of this subspace is finite. From the first half of the proof,
we have a lower estimate of the dimension of $J^{(k,r)}_M$.
We give an upper estimate of the dimension of the same space
by considering its dual space. This is a standard technique originated in
the paper by Feigin and Stoyanovsky \cite{FS}. Showing that these two estimates
are equal, we finish the proof of Theorem \ref{NEW}.

%%%%%%%%%%          Section koornwinder          %%%%%%%%%%

\section{Properties of the Koornwinder-Macdonald polynomials}
~
Let $n$ be the number of variables.
We denote by $W_n$ the group generated by permutations and sign flips
 ($W_n \cong {\mathfrak S}_n \ltimes (\mathbb{Z}_2)^n$).
We consider a $W_n$-symmetric Laurent polynomial ring
\bea
\bar{\Lambda}_n=\mathbb{C}[x_1^{\pm 1},\cdots ,x_n^{\pm 1}]^{W_n}.\label{BARL}
\eea

We denote by $\pi_n$ the set of partitions of length $n$, $\lambda=(\lambda_1,\cdots,\lambda_n)$.
We denote by $\widehat{m}_\lambda$ a monomial $W_n$-symmetric Laurent polynomial:
\bea
\widehat{m}_\lambda(x)&:=&\sum_{\nu\in W_n\lambda}\prod_ix_i^{\nu_i}.\nonumber
\eea

Let $\Lambda_n=\bar{\Lambda}_n\otimes\mathbb{C}(t,q,a,b,c,d)$.
The Koornwinder-Macdonald polynomial $P_\lambda(x)$ corresponding to $\lambda$
 is a simultaneous eigenfunction of the difference operators
 $\{D_r; 1\leq r\leq n\}$ (see \cite{vD}).
The corresponding eigenvalues $E^{(r)}_\lambda$ are of the form
\bea
E^{(r)}_\lambda:=u_\lambda(\widehat{m}_{1^{r}})+\sum_{0\leq s<r}a_{r,s} u_\lambda(\widehat{m}_{1^{s}}) \nonumber
\eea
where $u_\lambda$ is the one in (\ref{dual_param}) and $a_{r,s}\in$$\mathbb{C}[t^{\pm1},q^{\pm1},a^{\pm1},b^{\pm1},c^{\pm1},d^{\pm1},(a^*)^{\pm1}]$.

To be precise,
\bea
D_r&:=&\sum_{\substack{ J\subset\{1,\cdots,n\},0\leq|J|\leq r \\ \epsilon_j=\pm1,j\in J }} U_{J^c,r-|J|}(x)V_{\epsilon J,J^c}(x)T_{\epsilon J,q} \nonumber
\eea
\bea
V_{\epsilon J,K}(x)&:=&\prod_{j\in J} a^*
        \frac{1-ax_j^{\epsilon_j}}{1-x_j^{\epsilon_j}}
        \frac{1-bx_j^{\epsilon_j}}{1+x_j^{\epsilon_j}}
        \frac{1-cx_j^{\epsilon_j}}{1-q^{1/2}x_j^{\epsilon_j}}
        \frac{1-dx_j^{\epsilon_j}}{1+q^{1/2}x_j^{\epsilon_j}}\nonumber\\
&\times&\prod_{j,j'\in J,j<j'} t^{-1}
       \frac{1-tx_j^{\epsilon_j}x_{j'}^{\epsilon_{j'}}}{1-x_j^{\epsilon_j}x_{j'}^{\epsilon_{j'}}} 
       \frac{1-tqx_j^{\epsilon_j}x_{j'}^{\epsilon_{j'}}}{1-qx_j^{\epsilon_j}x_{j'}^{\epsilon_{j'}}}  \nonumber\\
&\times&\prod_{j\in J,k\in K} t^{-1}
       \frac{1-tx_j^{\epsilon_j}x_k}{1-x_j^{\epsilon_j}x_k} 
       \frac{1-tx_j^{\epsilon_j}x_k^{-1}}{1-x_j^{\epsilon_j}x_k^{-1}}  \nonumber\\
U_{K,p}(x)&:=&(-1)^p\sum_{\substack{L\subset K,|L|=p \\ \epsilon_l=\pm1,l\in L}}\prod_{l\in L} a^*
        \frac{1-ax_l^{\epsilon_l}}{1-x_l^{\epsilon_l}}
        \frac{1-bx_l^{\epsilon_l}}{1+x_l^{\epsilon_l}}
        \frac{1-cx_l^{\epsilon_l}}{1-q^{1/2}x_l^{\epsilon_l}}
        \frac{1-dx_l^{\epsilon_l}}{1+q^{1/2}x_l^{\epsilon_l}}\nonumber\\
&\times&\prod_{l,l'\in L,l<l'} t^{-1}
    \frac{1-tx_l^{\epsilon_l}x_{l'}^{\epsilon_{l'}}}{1-x_l^{\epsilon_l}x_{l'}^{\epsilon_{l'}}}
    \frac{1-tq^{-1}x_l^{-\epsilon_l}x_{l'}^{-\epsilon_{l'}}}{1-q^{-1}x_l^{-\epsilon_l}x_{l'}^{-\epsilon_{l'}}}   \nonumber\\
&\times&\prod_{l\in L,k\in K\backslash L} t^{-1}
    \frac{1-tx_l^{\epsilon_l}x_k}{1-x_l^{\epsilon_l}x_k}
    \frac{1-tx_l^{\epsilon_l}x_k^{-1}}{1-x_l^{\epsilon_l}x_k^{-1}}  \nonumber
\eea
and
\bea
%a_{r,s}&:=&(-1)^{r-s}\!\!\!\!\!\!\!\!\!\sum_{r \leq l_1 \leq \cdots \leq l_{r-s} \leq n}
% \!\!\!\!\!\!\!\!\!\!(t^{n-l_1}a^*+t^{-n+l_1}(a^*)^{-1})\cdots (t^{n-l_{r-s}}a^*+t^{-n+l_{r-s}}(a^*)^{-1}), \nonumber
a_{r,s}&:=&(-1)^{r-s}\sum_{r \leq l_1 \leq \cdots \leq l_{r-s} \leq n}
 \prod_{i=1}^{r-s}(t^{n-l_i}a^*+t^{-n+l_i}(a^*)^{-1}), \nonumber
\eea
where $T_{\epsilon J,q}:=\prod_{j\in J}T_{\epsilon_j j,q}$ , and
\bea
(T_{\pm j,q}f)(x_1,\cdots,x_n):=f(x_1,\cdots,x_{j-1},q^{\pm1}x_{j},x_{j+1},\cdots,x_n). \nonumber
\eea

For an indeterminate $X$, by taking the linear combination of $\{D_r\}$, we can define the operator $D(X)$
\bea
D(X):= \sum_{i=0}^{n} D'_iX^{n-i}, \nonumber
\eea
where $\{D'_i;0\leq i \leq n\}$ are defined inductively as follows
\bea
D'_0  &=&1 \nonumber\\
D'_{i}&=&D_{i}-\sum_{j<i} a_{i,j}D'_j . \nonumber
\eea
Then the eigenvalue $E_\lambda(X)$ of the operator $D(X)$ is given by
\bea
D(X)P_\lambda&=&E_\lambda(X)P_\lambda \nonumber\\
E_\lambda(X) &:=&\prod_{i=1}^{n}
 (X+t^{n-i}q^{\lambda_i}a^*+t^{-n+i}q^{-\lambda_i}(a^*)^{-1}).\nonumber
\eea

We use the dominance ordering $\lambda>\mu$ for partitions $\lambda$ and $\mu$.
We have

\begin{lem}\label{lem:welldef}
Let $c_{\lambda\mu}$ be
\bea
P_{\lambda}=:\widehat{m}_\lambda+\sum_{\mu<\lambda}c_{\lambda\mu}\widehat{m}_{\mu}.
\nonumber
\eea
If there does not exist $\nu<\lambda$ such that $E_\lambda(X)=E_\nu(X)$
 at a certain specialization of parameters,
 then for any $\mu<\lambda$, $c_{\lambda\mu}$ has no pole at the same specialization.
\end{lem}

\begin{proof}
It is clear from the defining equality of $P_{\lambda}$
\bea
P_{\lambda}:=\left(\prod_{\mu<\lambda}\frac{D(X)-E_\mu(X)}{E_\lambda(X)-E_\mu(X)}\right)m_\lambda. \nonumber
\eea
\end{proof}

\begin{rem}\normalfont
In order to distinguish eigenvalues in the specialization (\ref{resonance}),
 the second order operator $D_1$ alone is not enough.
The eigenvalue $E_\lambda^{(1)}$ of $D_1$ is given by
\[E_\lambda^{(1)}=\sum_{i=1}^{n}
 (t^{n-i}q^{\lambda_i}a^*+t^{-n+i}q^{-\lambda_i}(a^*)^{-1}).\]
For example, let $n=4$, $k=3$, and $r=3$. Then $t^2q=-1$.
Hence $E_{(3,3,2,0)}^{(1)}=E_{(4,3,3,0)}^{(1)}$,
 although $E_{(3,3,2,0)}(X) \neq E_{(4,3,3,0)}(X)$.
This is why we use the operator $D(X)$.
\end{rem}

We define a dual Koornwinder-Macdonald polynomial $P^*_\lambda$
 by dual parameters $a^*,b^*,c^*,d^*$ given in (\ref{dual_param}).

{}In \cite{Sa}, we have the following relation.

\begin{prop}[duality]
For all $\lambda,\mu \in \pi_n$,
 the Koornwinder-Macdonald polynomial $P_\lambda$
 and the dual Koornwinder-Macdonald polynomials $P^*_\mu \in \Lambda_n$ satisfy the following duality relation:
\bea
\frac{u^*_\mu(P_\lambda)}{u^*_0(P_\lambda)} = \frac{u_\lambda(P^*_\mu)}{u_0(P^*_\mu)}.\label{duality}
\eea
Here, the definition of $u^*_\mu$ and $u_\lambda$ are those in $(\ref{special_u})$.
\end{prop}

In \cite{vD}, it is shown that the duality relation (\ref{duality})
 implies the following evaluation formula.
\begin{prop}
\bea
&&u^*_0(P_\lambda)=P_\lambda^{sum} \times P_\lambda^{diff} \times P_\lambda^{single},\nonumber\\
&&P_\lambda^{sum}:= \prod_{i<j} t^{-(\lambda_i+\lambda_j)/2}
 \frac{(t^{2n+1-i-j}(a^*)^2;q)_{\lambda_i+\lambda_j}}
      {(t^{2n-i-j}(a^*)^2;q)_{\lambda_i+\lambda_j}} \label{sum} ,\\
&&P_\lambda^{diff}:= \prod_{i<j} t^{-(\lambda_i-\lambda_j)/2}
 \frac{(t^{j-i+1};q)_{\lambda_i-\lambda_j}}
      {(t^{j-i};q)_{\lambda_i-\lambda_j}} \label{diff},\\
&&P_\lambda^{single}:= \prod_{i} a^{-\lambda_i} \frac{(t^{n-i}(a^*)^2,t^{n-i}a^*b^*,t^{n-i}a^*c^*,t^{n-i}a^*d^*;q)_{\lambda_i}}
      {(t^{n-i}a^*,-t^{n-i}a^*,t^{n-i}a^*q^{1/2},-t^{n-i}a^*q^{1/2};q)_{\lambda_i}}.\label{single}
\eea
Here, $(a;q)_l:=\prod_{i=0}^{l-1}(1-aq^i)$ and $(a_1,a_2,\cdots,a_p;q)_l:=\prod_{i=1}^{p}(a_i;q)_l$.
\end{prop}

\begin{rem}\label{rem:norm}\normalfont
Note that in (\ref{diff}), there appear only factors of the form $(1-t^xq^y)$, $x,y\in\mathbb{Z}_{\geq0}$.
In (\ref{sum}), there appear only factors of the form $(1-t^xq^y(a^*)^2)$, $x,y\in\mathbb{Z}_{\geq0}$.
In (\ref{single}), there appear only factors of the form $(1-t^xq^y(a^*)^2)$,
 $(1-t^xq^ya^*b^*)$, $(1-t^xq^ya^*c^*)$, $(1-t^xq^ya^*d^*)$, $x,y\in\mathbb{Z}_{\geq0}$.
\end{rem}

%%%%%%%%%%          Section lower          %%%%%%%%%%

\section{The space $I^{(k,r)}_M$ and $J^{(k,r)}_M$}
~
In this section,
 we describe zero conditions and construct symmetric Laurent polynomials satisfying the zero conditions.

First, we describe a specialization of the parameters.
Let $k,r$ be integers such that $1\leq k \leq n-1$ and $r\geq 2$.
Let $m$ be the greatest common divisor of $(k+1)$ and $(r-1)$.
Let $\omega$ be a primitive $m$-th root of unity.
Let $\omega_1\in\mathbb{C}$ be such that $\omega_1^{(r-1)/m}=\omega$.

\begin{defn}\normalfont
For an indeterminate $u$, we consider the specialization of $t$ and $q$:
\bea
t=u^{(r-1)/m}, q=\omega_1u^{-(k+1)/m}.\label{resonance}
\eea

Then for integers $x,y\in\mathbb{Z}$,
 $t^xq^y=1$ if and only if $x=(k+1)l$, $y=(r-1)l$ for some $l\in\mathbb{Z}$.
Moreover, the multiplicity of $(t^{(k+1)/m}q^{(r-1)/m}-\omega)$ in $(t^{(k+1)l}q^{(r-1)l}-1)$ is $1$.
\end{defn}

We define the subject of our study.
We denote by $\Lambda_n'$ the corresponding space $\Lambda_n':=\bar{\Lambda}_n\otimes\mathbb{C}(u,a,b,c,d)$.

\begin{defn}\normalfont
A sequence $(s_1,\cdots,s_{k+1})$ ($s_1,\cdots, s_{k+1}\in \mathbb{Z}_{\geq0}$) is called a $wheel$ $sequence$ if $s_1+\cdots+s_{k+1}=r-1$.
For $f \in \Lambda_n'$, we consider the following $wheel\ condition$:
\begin{eqnarray}
\begin{array}{l}
f=0, \mbox{\quad if } x_{i+1}=tq^{s_i}x_i \qquad (1\leq i\leq k) \\
\mbox{for all wheel sequences $(s_1,\cdots,s_{k+1})$}.\label{WC}
\end{array}
\end{eqnarray}
We consider the subspace $J^{(k,r)}\subseteq\Lambda_n'$
\bea
J^{(k,r)}:=\{ f\in\Lambda_n' ; \mbox{$f$ satisfies (\ref{WC})} \}.
\eea

We denote by $\Lambda_{n,M}'$ the subspace consisting of $f\in\Lambda_n'$ such that the degree of $f$ in each $x_i$ is less than $M$.
We set $J^{(k,r)}_M:= J^{(k,r)}\cap \Lambda_{n,M}'$.
\end{defn}

\begin{rem}\label{rem:wheel}\normalfont
For any partition $\mu\in\pi_n$,
 $u^*_\mu(x_1)/u^*_\mu(x_{k+1})=t^{k}q^{\mu_1-\mu_{k+1}}$.
Hence the condition $\mu_1-\mu_{k+1}\leq r-1$ corresponds to
 the existence of the wheel sequence: $s_{k+1}=r-1-(\mu_1-\mu_{k+1})\geq 0$.
The wheel conditions for $f(x)\in\Lambda_n'$ correspond to
 $u^*_\mu(f)=0$ at the specialization $(\ref{resonance})$
 for any partition $\mu\in\pi_n$ such that $\mu_1-\mu_{k+1}\leq r-1$.
\end{rem}

\begin{rem}\normalfont
In \cite{DS}, it is proved that
 $u^*_\mu(P_\lambda)=0$ for $\mu_1\leq N$ and $\lambda_1> N$
 under the specialization of parameters $t^{n-1}abq^N=1$.
As it were, this is the zero condition on the finite set
 $\{(u^*_\mu(x_1),\cdots,u^*_\mu(x_n));\mu\in\pi_n,\mu_1\leq N\}$.
On the other hand, in this paper, the wheel condition is on the infinite set,
 because it determines at most the ratio of variables.
\end{rem}

For $f(t,q,a,b,c,d)\in\mathbb{C}[t,q,a,b,c,d]$, we use a specialization mapping $\varphi$
\bea
\varphi \mbox{ : } \mathbb{C}[t,q,a,b,c,d] &\longrightarrow& 
\mathbb{C}(u,a,b,c,d)\nonumber \\
f(t,q,a,b,c,d) &\mapsto& f(u^{(r-1)/m},\omega_1u^{-(k+1)/m},a,b,c,d),\nonumber
\eea
and we extend $\varphi$ to those elements of the field
 $\mathbb{C}(t,q,a,b,c,d)$ for which the specialized denominator does not vanish.

\begin{lem}\label{lem:welldef2}
If $\lambda\in\pi_n$ satisfies
\bea
\lambda_i-\lambda_{i+k+1}> 2\left[\frac{n}{k+1}\right](r-1) \quad\mbox{for $1\leq i \leq n-k-1$},\nonumber
\eea
then $P_\lambda$ and $P^*_\lambda$ have no pole at the specialization $(\ref{resonance})$.
\end{lem}

\begin{proof}
First, we discuss on $P_\lambda$.
Suppose that there exists $\mu$ such that $\varphi(E_\mu(X))=\varphi(E_\lambda(X))$, that is
\bea
&&\{ \varphi( t^{n-i}q^{\mu_i}a^*+t^{-n+i}q^{-\mu_i}(a^*)^{-1} ) ; 1\leq i \leq n \} \nonumber\\
&&=\{ \varphi( t^{n-i}q^{\lambda_i}a^*+t^{-n+i}q^{-\lambda_i}(a^*)^{-1} ) ; 1\leq i \leq n \}.
\nonumber
\eea
Since $u$ and $a^*$ are generic, it must be satisfied that
\bea
&&\{ \varphi( t^{n-(k+1)l-i}q^{\mu_{(k+1)l+i}} ) \mbox{ ; $l\in\mathbb{Z}_{\geq0}$ and $1 \leq (k+1)l+i \leq n$} \}\nonumber\\
&&=\{ \varphi( t^{n-(k+1)l-i}q^{\lambda_{(k+1)l+i}}) \mbox{ ; $l\in\mathbb{Z}_{\geq0}$ and $1 \leq (k+1)l+i \leq n$} \}\nonumber
\eea
for $1\leq i \leq k+1$.
Hence
\bea
&&\{ (r-1)l+\mu_{(k+1)l+i} \mbox{ ; $l\in\mathbb{Z}_{\geq0}$ and $1 \leq (k+1)l+i \leq n$} \}\nonumber\\
&&= \{ (r-1)l+\lambda_{(k+1)l+i} \mbox{ ; $l\in\mathbb{Z}_{\geq0}$ and $1 \leq (k+1)l+i \leq n$} \}\nonumber
\eea
for $1\leq i \leq k+1$.

Then for any $1\leq i \leq k+1$,
 there exists $l_i \geq 0$ such that $(r-1)l_i+\mu_{(k+1)l_i+i}=\lambda_i$
 and there exists $l'_i \geq 0$ such that $(r-1)l'_i+\lambda_{(k+1)l'_i+i}=\mu_i$.
If $l'_i \neq 0$, then by the hypothesis,
\bea
\mu_i-\mu_{(k+1)l_i+i} &=& (r-1)l'_i+\lambda_{(k+1)l'_i+i}-\lambda_i+(r-1)l_i \nonumber\\
&<& (r-1)(l'_i+l_i) -2\left[\frac{n}{k+1}\right](r-1)l'_i \nonumber\\
&\leq& 2\left[\frac{n}{k+1}\right](r-1)(1-l'_i) \nonumber\\
&\leq& 0. \nonumber
\eea
Hence $l'_i$ must be equal to $0$, namely $\lambda_i=\mu_i$.
Inductively, we have $\lambda_{(k+1)l+i}=\mu_{(k+1)l+i}$ for all $l\geq 0$.
It follows that $\lambda=\mu$.
Therefore from Lemma \ref{lem:welldef}, $P_\mu$ has no pole at the specialization (\ref{resonance}).

For $P^*_\lambda$, its eigenvalue $E^*_\lambda(X)$ is given by replacing $a^*$
 with $a$ in $E_\lambda(X)$.
Hence we can similarly show that $P^*_\lambda$ has no pole at the specialization (\ref{resonance}).
\end{proof}

We are going to construct a basis of $J^{(k,r)}_M$.

\begin{defn}\normalfont
$\lambda \in \pi_n$ is called $(k,r,n)$-$admissible$ if
\begin{eqnarray}
\lambda_i-\lambda_{i+k}\geq r \qquad(1\leq \forall i\leq n-k).
\end{eqnarray}
\end{defn}

Our main result is

\begin{thm}\label{thm:main}
For any $(k,r,n)$-admissible $\lambda$,
 Koornwinder-Macdonald polynomial $P_\lambda$ has no pole
 at the specialization $(\ref{resonance})$.
Moreover, for any positive integer $M$, we have
\bea
I^{(k,r)}_M=J^{(k,r)}_M. 
\nonumber
\eea
Here, we define a subspace $I^{(k,r)}$ of $\Lambda_n'$
\bea
I^{(k,r)}:=\operatorname{span}_{\mathbb{C}(u,b,c,d)}\{ \varphi(P_\lambda); \lambda \mbox{ is $(k,r,n)$-admissible } \},
\nonumber
\eea
and we set
\bea
I^{(k,r)}_M:=\operatorname{span}_{\mathbb{C}(u,b,c,d)}\{ \varphi(P_\lambda); \lambda \mbox{ is $(k,r,n)$-admissible and $\lambda_1\leq M$} \}.
\nonumber
\eea

\end{thm}

First, we prepare some propositions and lemmas.

\begin{defn}\normalfont
For $p\in\mathbb{C}(t,q,a,b,c,d)$, we denote by $Z(p)\in\mathbb{Z}$ the multiplicity of $(t^{(k+1)/m}q^{(r-1)/m}-\omega)$ in $p$. That is,
\bea
p = (t^{k+1}q^{r-1}-1)^{Z(p)}p',\nonumber
\eea
where the factor $p'\in\mathbb{C}(t,q,a,b,c,d)$ has neither pole nor zero at $(\ref{resonance})$.
\end{defn}

\begin{prop}\label{prop:zeros and poles}
For any partition $\lambda\in\pi_n$, we have
\bea
Z(u^*_0(P_\lambda))=Z(u_0(P^*_\lambda))
 &=& \sharp\{(i,l)\in\mathbb{Z}_{> 0}^2; \lambda_i-\lambda_{i+(k+1)l-1}\geq (r-1)l+1 \}\nonumber\\
 &&- \sharp\{(i,l)\in\mathbb{Z}_{> 0}^2; \lambda_i-\lambda_{i+(k+1)l}\geq (r-1)l+1 \}.\nonumber
\eea
\end{prop}

\begin{proof}
Recall Remark \ref{rem:norm}.
The factor $P_\lambda^{diff}$ has the factors of the form $(1-t^xq^y)$ ($x,y\in\mathbb{Z}_{\geq 0}$).

If $j-i+1=(k+1)l$ and $\lambda_i-\lambda_j\geq (r-1)l+1$,
 then $u^*_0(P_\lambda)$ has the factor $(1-t^{(k+1)l}q^{(r-1)l})$ in the numerator of $P_\lambda^{diff}$.
If $j-i=(k+1)l$ and $\lambda_i-\lambda_j\geq (r-1)l+1$,
 then $u^*_0(P_\lambda)$ has the factor $(1-t^{(k+1)l}q^{(r-1)l})$ in the denominator of $P_\lambda^{diff}$.
Otherwise, there does not exist the factor $(1-t^{(k+1)l}q^{(r-1)l})$ in $P_\lambda^{diff}$.

On the other hand,
 $P_\lambda^{sum}$ and $P_\lambda^{single}$ have neither pole nor zero at the specialization $(\ref{resonance})$.

$u_0(P^*_\lambda)$ is given by replacing parameters with dual ones in $u^*_0(P_\lambda)$.
Since the specialization $(\ref{resonance})$ is invariant under
 acting $*$ on parameters, we have $Z(u_0(P^*_\lambda))=Z(u^*_0(P_\lambda))$.
\end{proof}

\begin{cor}\label{cor:adm}
For any $(k,r,n)$-admissible $\lambda$, we have $Z(u^*_0(P_\lambda))=[\frac{n}{k+1}]$.
\end{cor}

\begin{proof}
Since $\lambda_i-\lambda_{i+k}\geq r$,
\bea
Z(u^*_0(P_\lambda))
 &=& \sharp\{(i,l)\in\mathbb{Z}_{> 0}^2; i+(k+1)l-1\leq n \} \nonumber \\
 &&- \sharp\{(i,l)\in\mathbb{Z}_{> 0}^2; i+(k+1)l  \leq n \} \nonumber \\
 &=& \sum_{l\geq1} \operatorname{max}\{(n-(k+1)l+1),1\}
 - \sum_{l\geq1} \operatorname{max}\{(n-(k+1)l),1\} \nonumber \\
 &=&\left[\frac{n}{k+1}\right] \nonumber
\eea
\end{proof}

\begin{rem}\label{rem:to prove}\normalfont
%Let $m$ be the greatest common divisor of $(k+1)$ and $(r-1)$.
%Let $\omega$ be a primitive $m$-th root of unity and $F(X)$ be
%the minimal polynomial of $\omega$ over $\mathbb{C}$.
%Set $f:=F(t^{(k+1)/m}q^{(r-1)/m})\in \mathbb{C}[q,t,b,c,d]$.
%Let $I$ be the ideal of the polynomial ring $\mathbb{C}[q,t,b,c,d]$
%generated by $f$.
%
%To show that $f\in\Lambda_n$ has no pole at the specialization
%(\ref{resonance}),
% it is sufficient to check $p \cdot f = 0$ at the specialization
%(\ref{resonance}) for all $p\in I$.

For $g\in\Lambda_n$, we take an integer $N$ such that the degree
of $g$ in each variable $x_i$ is less than $N/2$.
Then to prove that $g = 0$ (respectively, $g$ has no pole) at the specialization (\ref{resonance}),
 it is sufficient to show that there exist $n$ subsets
$C_1,\cdots,C_n \subseteq \mathbb{C}(a,b,c,d)[q^{\pm1},t^{\pm1}]$,
 which satisfy the following two conditions:
\par
(i) for each $i$, $\sharp(\varphi(C_i))\geq N$ in $\mathbb{C}(u,a,b,c,d)$;
\par
(ii) for all choices of $c_i\in C_i$ , $Z(g(c_1,\cdots,c_n))>0$ (resp. $\geq 0$).
\end{rem}

Motivated by the observation above, we define certain sets of partitions.

\begin{defn}\normalfont
A partition $\eta$ is called thick if $\eta_i \gg \eta_{i+1} \gg 0$ for all $i$.
For a thick partition $\eta\in\pi_n$, a set of $N^n$ partitions is defined by
$\pi_{\eta,N}:=\{\mu\in\pi_n ; \mu_i=\eta_i+d_i$ for all $i$ where $0\leq d_i\leq N-1\}$. 

For a thick partition $\eta\in\pi_{n-k}$, we define
$\pi_{\eta,N}':=\{\mu\in\pi_n ;\mu_1-\mu_{k+1}<r, \mu_i=\eta_{i-k}+d_{i-k}$ for $k+1\leq i \leq n$ where $0\leq d_i\leq N-1\}$. 
\end{defn}

When we use these sets $\pi_{\eta,N}$ and $\pi_{\eta,N}'$,
 we choose a sufficiently large $N$ such that $N \gg M$
 and any thick partition $\eta$ such that
 $\eta_i-\eta_{i+1}\gg \operatorname{max}(M, 2[\frac{n}{k+1}](r-1))$,
 $\eta_i\gg \operatorname{max}(M, 2[\frac{n}{k+1}](r-1))$.
We do not specify $N$ and $\eta$ in the below.

\begin{lem}\label{lem:thickwheel}
For $\mu\in\pi_{\eta,N}$ or $\mu\in\pi_{\eta,N}'$,
 $P_\mu$ and $P^*_\mu$ have no pole at the specialization $(\ref{resonance})$.
Moreover
\bea
Z(u^*_0(P_\mu)) = Z(u_0(P^*_\mu)) = 
 \begin{cases}
 [\frac{n}{k+1}] &\mbox{ if } \mu\in\pi_{\eta,N} \ (\eta\in\pi_n), \\[0pt]
 [\frac{n}{k+1}]-1 &\mbox{ if } \mu\in\pi_{\eta,N}' \ (\eta\in\pi_{n-k}).
 \end{cases}
\nonumber
\eea
\end{lem}

\begin{proof}
If $\mu$ is an element of $\pi_{\eta,N}$ or $\pi_{\eta,N}'$,
 then $\mu_i\gg\mu_{i+k+1}$ for $1\leq i \leq n-k-1$.
Hence from Lemma \ref{lem:welldef2}, we see $P_\mu$ and $P^*_\mu$
 have no pole at $(\ref{resonance})$.

If $\mu\in\pi_{\eta,N}$ , then for each $1 \leq l \leq [\frac{n}{k+1}]$,
 $\mu_i\gg\mu_{i+(k+1)l-1}$ ($1\leq i\leq n-(k+1)l+1$) and $\mu_i\gg\mu_{i+(k+1)l}$ ($1\leq i\leq n-(k+1)l$).
Hence from Proposition \ref{prop:zeros and poles},
 $Z(u^*_0(P_\mu)) = Z(u_0(P^*_\mu)) = [\frac{n}{k+1}]$.

If $\mu\in\pi_{\eta,N}'$, then $\mu_1-\mu_{k+1}\leq r-1$.
Hence from Proposition \ref{prop:zeros and poles},
 $(i,l)=(1,1)$ is the only different situation from the case $\mu\in\pi_{\eta,N}$.
Therefore $Z(u^*_0(P_\mu)) = Z(u_0(P^*_\mu)) = [\frac{n}{k+1}]-1$.
\end{proof}

Now we are ready to prove a part of Theorem \ref{thm:main}.

\begin{thm}\label{thm:generic}
For any $(k,r,n)$-admissible $\lambda$,
 Koornwinder-Macdonald polynomial $P_\lambda$ has no pole at the specialization $(\ref{resonance})$
 and $\varphi(P_\lambda)$ satisfies the wheel condition $(\ref{WC})$.
\end{thm}

\begin{proof}
Since $\lambda$ is $(k,r,n)$-admissible,
 $Z(u^*_0(P_\lambda)) = [\frac{n}{k+1}]$ from Corollary \ref{cor:adm}.

Let $N\gg |\lambda|$ and let $\mu\in\pi_{\eta,N}$ where $\eta\in\pi_n$.
Then from Lemma \ref{lem:thickwheel},
 $P^*_\mu$ has no pole at the specialization (\ref{resonance}) and
 $Z(u_0(P^*_\mu))=[\frac{n}{k+1}]$.
From the duality relation (\ref{duality}),
\bea
u^*_\mu(P_\lambda)=\frac{u_\lambda(P^*_\mu)}{u_0(P^*_\mu)}u^*_0(P_\lambda).
\nonumber
\eea
Therefore, $Z(u^*_\mu(P_\lambda))\geq0$.

Since this holds for all $\mu\in\pi_{\eta,N}$,
 from Remark \ref{rem:to prove},
 we see that $P_\lambda$ has no pole at the specialization (\ref{resonance}).

Let $\mu\in\pi_{\eta,N}'$ ($\eta\in\pi_{n-k}$).
Then from Lemma \ref{lem:thickwheel},
 $P^*_\mu$ has no pole at the specialization (\ref{resonance}) and $Z(u_0(P^*_\mu))=[\frac{n}{k+1}]-1$.
From the duality relation (\ref{duality}),
 through the same argument as the above, $Z(u^*_\mu(P_\lambda))\geq 1$.

We have shown $u^*_\mu(P_\lambda)=0$ at the specialization (\ref{resonance}) for all $\mu\in\pi_{\eta,N}'$.
Hence from Remark \ref{rem:wheel} and Remark \ref{rem:to prove},
 we conclude that $\varphi(P_\lambda)$ satisfies the wheel condition (\ref{WC}).

\end{proof}

\begin{cor}\label{cor:lower}
The space $I^{(k,r)}$ and $I^{(k,r)}_M$ are well-defined for any positive integer $M$,
 and we have
$J^{(k,r)}_M \supseteq I^{(k,r)}_M$.
\end{cor}

%%%%%%%%%%          section estimate          %%%%%%%%%%

\section{Estimate of $\dim J^{(k,r)}_M$}
We have already constructed the polynomials satisfying the zero conditions.
In this section,
 we show that $J^{(k,r)}_M = I^{(k,r)}_M$ by giving an upper estimate of the dimension of $J^{(k,r)}_M$.

Fix $g_0'$, $g_1'$, $g_2'$, $g_3' \gg 1$.
We take the limit $t\rightarrow 1$, $q\rightarrow \tau$,
 $a\rightarrow \tau^{g_0'}$, $b\rightarrow -\tau^{g_1'}$,
 $c\rightarrow \tau^{g_2'+1/2}$, $d\rightarrow -\tau^{g_3'+1/2}$,
 where $\tau$ is a primitive $(r-1)$-th root of unity.
In this limit the wheel condition (\ref{WC}) reduces to
\begin{eqnarray}
f=0&&\mbox{ if } x_{i}=\tau^{p_{i}}x_0 \qquad (1\leq i\leq k+1)\label{WC0}
\end{eqnarray}
for all $p_1,\cdots,p_{k+1}\in \mathbb{Z}$ and $x_0\in\mathbb{C}$.
We denote by $\bar{J}^{(k,r)}\subseteq\bar{\Lambda}_n$
the space of $(BC)_n$-symmetric polynomials satisfying (\ref{WC0}).
We define
\bea
\bar{J}^{(k,r)}_M=\{f\in\bar{J}^{(k,r)} ; \operatorname{deg}_{x_1}f\leq M \}. \nonumber
\eea
Note that $\dim_{\mathbb{C}(u,b,c,d)} J^{(k,r)}_M \leq \dim_{\mathbb{C}}
\bar{J}^{(k,r)}_M$.

We consider the commutative ring $R_M:=\mathbb{C}[e_0,e_1,e_2,\cdots,e_M]$ for indeterminates $\{e_i\}$.
We count the weight of $e_i$ as $1$ and the degree of $e_i$ as $i$.
We set $e_\lambda:=\prod_{i=1}^{n}e_{\lambda_i}$ for $\lambda \in \pi_n$.
We denote by $R_{M,n}\subseteq R_{M}$ the space spanned by the monomials $e_\lambda$ such that $\lambda\in\pi_n$ and $\lambda_1\leq M$.

We use the dual language (see \cite{FJMM}).
There is a nondegenerate coupling:
\bea
R_{M,n}\times\bar{\Lambda}_{n,M}\rightarrow \mathbb{C};\\
\langle e_\lambda , \widehat{m}_\mu \rangle=\delta_{\lambda,\mu}.
\nonumber
\eea

We introduce an abelian current
\bea
e(z)&:=&\sum_{i=1}^{M} e_i(z^i+z^{-i}) + e_0.
\nonumber
\eea
It satisfies
\bea
\langle e(z_1)e(z_2) \cdots e(z_n),f \rangle = f(z_1,z_2,\cdots,z_n) \quad\mbox{for $f\in\bar{\Lambda}_{n,M}$}.\nonumber
\eea
Then for any $f\in\bar{J}^{(k,r)}_M$, we have
\bea
\langle e(\tau^{p_1} z) \cdots e(\tau^{p_{k+1}} z)e(z_{k+2})\cdots e(z_n),f \rangle=0
 \quad\mbox{for all}\ (p_1,\cdots,p_{k+1})\in \mathbb{Z}^{k+1}. \nonumber
\eea
Hence the space
\bea
\begin{array}{ll}
 &\operatorname{span}_\mathbb{C}\{e(\tau^{p_1} z) \cdots e(\tau^{p_{k+1}} z)e(z_{k+2})\cdots e(z_n) \\
 &\qquad\qquad\qquad;\ z,z_{k+2},\cdots,z_n\in\mathbb{C},p_1,\cdots,p_{k+1}\in\mathbb{Z}\}
\end{array}\label{ortho}
\eea
is the orthogonal complement of $\bar{J}^{(k,r)}_M$ with respect to the coupling $\langle,\rangle$.
For $p=(p_1,\cdots,p_{k+1})\in \mathbb{Z}^{k+1}$, let $r_d^p$ be the coefficient of $z^d$ in
\bea
e(\tau^{p_1} z)\cdots e(\tau^{p_{k+1}} z)=\sum_d r_d^pz^d.
\nonumber
\eea
By the symmetry of exchanging $z\leftrightarrow z^{-1}$ in the current $e(z)$, we have $r_d^p = r_{-d}^p$.
We denote by ${\cal J}_M$ the ideal of $R_M$ generated by the elements $r_d^p$.
Set ${\cal J}_{M,n}:={\cal J}_M\cap R_{M,n}$.
Then the space (\ref{ortho}) coincides with ${\cal J}_{M,n}$.
Since $\dim R_{M,n}/{\cal J}_{M,n} = \dim \bar{J}^{(k,r)}_M$,
 the condition (\ref{WC0}) is equivalent to the relations in the quotient space
\bea
r_d^{p}=0 \quad\mbox{ for all $p=(p_1,\cdots,p_{k+1})\in\mathbb{Z}^{k+1}$ and $d\geq 0$}.\nonumber
\eea

\begin{prop}
The image of the set $\{e_\lambda;\lambda\in \pi_n$ is $(k,r,n)$-admissible,$\lambda_1\leq M\}$ spans the quotient space $R_{M,n}/{\cal J}_{M,n}$.
\end{prop}

\begin{proof}
We introduce a total ordering for partitions and monomials.
For two partitions $\lambda$ and $\mu$ such that $|\lambda| > |\mu|$,
 we define $\lambda \succ \mu$.
For two partitions $\lambda$ and $\mu$ such that $|\lambda| = |\mu|$,
 we define $\lambda \succ \mu$ if $\lambda_1>\mu_1$ or $\lambda_1=\mu_1, \lambda_2>\mu_2$ or $\lambda_1=\mu_1, \lambda_2=\mu_2, \lambda_3>\mu_3$ or $\cdots$.
We define for the corresponding monomials $e_\lambda$ and $e_\mu$,
 $e_\lambda \succ e_\mu$.

Let us calculate $r_d^p$.
\bea
e(\tau^{p_1} z)\cdots e(\tau^{p_{k+1}} z)&=& \prod_{j=1}^{k+1}\sum_{i_j=-M}^{M} e_{|i_j|}(\tau^{p_j}z)^{i_j}\nonumber\\
&=& \sum_{d\in\mathbb{Z}}z^d \left( \sum_{\substack{i_1+\cdots+i_{k+1}=|d| \\ i_j\geq0}}\prod_{j=1}^{k+1}e_{i_j}\tau^{p_j i_j}+\sum_{\substack{\lambda\in\pi_{k+1} \\ |\lambda|>|d|}}c_{\lambda,|d|}e_\lambda \right). \nonumber
\eea
Hence, for any nonnegative integer $d$,
\bea
r_d^p=\sum_{\substack{i_1+\cdots+i_{k+1}=d \\ i_j\geq0}}\prod_{j=1}^{k+1}e_{i_j}\tau^{p_j i_j}+\sum_{\substack{\lambda\in\pi_{k+1} \\ |\lambda|>d}}c_{\lambda,d}e_\lambda. \nonumber
\eea
We define $R_{M,k+1}^{(d)}$ by
\bea
R_{M,k+1}^{(d)}:= \bigoplus_{\lambda\in\pi_{k+1}, |\lambda|\geq d, \lambda_1\leq M } \mathbb{C}e_\lambda,\nonumber
\eea
and we consider a quotient space
\bea
R_{M,k+1}^{(d)}/(R_{M,k+1}^{(d+1)}+\sum_{p\in\mathbb{Z}^{k+1}}\mathbb{C}r_d^p).\nonumber
\eea
In this space, 
\bea
0&=&r_d^p=\sum_{\substack{i_1+\cdots+i_{k+1}=d \\ i_j\geq0}}\prod_{j=1}^{k+1}e_{i_j}\tau^{p_j i_j} \nonumber\\
&=& \sum_{\substack{\nu\in\mathbb{Z}^{k+1}_{\geq0} \\ \nu_j\leq r-2}} \tau^{p_1\nu_1+\cdots+p_{k+1}\nu_{k+1}}\left(\sum_{\substack{\mu\in\mathbb{Z}^{k+1}_{\geq0}, \sum_j \mu_j=d \\ \mu_j=\nu_j+(r-1)\kappa_j , \kappa_j\in\mathbb{Z}_{\geq0}}}\prod_j e_{\mu_j} \right) \nonumber\\
&=& \sum_{\substack{\lambda\in\pi_{k+1} \\ \lambda_1\leq r-2}}\left(\sum_{\nu\in\mathfrak{S}_{k+1}\lambda} \tau^{p_1\nu_1+\cdots+p_{k+1}\nu_{k+1}}\right)\left(\sum_{\substack{\mu\in\mathbb{Z}^{k+1}_{\geq0}, \sum_j \mu_j=d \\ \mu_j=\lambda_j+(r-1)\kappa_j , \kappa_j\in\mathbb{Z}_{\geq0}}}\prod_j e_{\mu_j} \right)\nonumber
\eea
We set $\pi_{k+1,d}:=\{\lambda\in\pi_{k+1}; |\lambda|=d \}$.
For a sequence of nonnegative integers $m := (m_0,\cdots,m_{r-2})$ such that $\sum m_i = k+1$, we define a subset $\pi_{k+1,d}(m)$ by
\bea
&&\pi_{k+1,d}(m):=\{\mu \in \pi_{k+1,d} \ ; \nonumber\\
&&\ \qquad\qquad\qquad\sharp\{i;\mu_i\equiv a \ mod\ (r-1)\}=m_a \mbox{ for every }0\leq a \leq r-2\}.
\nonumber
\eea
We denote by $m_i^{(\lambda)}$ the multiplicity of $i$ in $\lambda$.
Define $m(\lambda):=(m_0^{(\lambda)},\cdots,m_{r-2}^{(\lambda)})$.
Then,
\bea
r_d^p&=&\sum_{\substack{\lambda\in\pi_{k+1} \\ \lambda_1\leq r-2}} \left( \sum_{\mu\in\pi_{k+1,d}\left(m(\lambda)\right)} c_{\lambda,\mu} e_\mu \right) \left( \sum_{\nu\in\mathfrak{S}_{k+1}\lambda} \tau^{p_1\nu_1+\cdots+p_{k+1}\nu_{k+1}} \right). \nonumber\\
&=&\sum_{\substack{\lambda\in\pi_{k+1} \\ \lambda_1\leq r-2}} \left( \sum_{\mu\in\pi_{k+1,d}\left(m(\lambda)\right)} c_{\lambda,\mu} e_\mu \right) m_\lambda(\tau^{p_1},\cdots,\tau^{p_{k+1}}) . \nonumber
\eea
Here, $c_{\lambda,\mu}=\prod_i m_i^{(\lambda)}!/\prod_i m_i^{(\mu)}!$
and $m_\lambda$ is the monomial $\mathfrak{S}_{k+1}$-symmetric polynomial (not Laurent).

Since $\lambda_1\leq r-2$, the degree of 
\bea
\sum_{\substack{\lambda\in\pi_{k+1} \\ \lambda_1\leq r-2}} \left( \sum_{\mu\in\pi_{k+1,d}\left(m(\lambda)\right)} c_{\lambda,\mu} e_\mu \right) m_\lambda(x_1,\cdots,x_{k+1}) \qquad \label{rel2}
\eea
in each variable $x_i$ is less than $r-2$.
On the other hand, we can choose the values of $x_i$ from $\tau^0,\tau^1,\cdots,\tau^{r-2}$ independently.
Hence the expression (\ref{rel2}) is identically zero in the quotient space
$R_{M,k+1}^{(d)}/(R_{M,k+1}^{(d+1)}+\sum_p \mathbb{C} r_d^p)$.
Since monomial symmetric polynomials are linearly independent, it follows that
\bea
\sum_{\mu\in\pi_{k+1,d}\left(m(\lambda)\right)} c_{\lambda,\mu} e_\mu=0 \nonumber
\eea
in $R_{M,k+1}^{(d)}/(R_{M,k+1}^{(d+1)}+\sum_p \mathbb{C} r_d^p)$.
Note that ${\cal J}_{M,k+1}=\sum_{d=0}^{M^{k+1}} \sum_p \mathbb{C} r_d^p$.
Therefore in $R_{M,k+1}/{\cal J}_{M,k+1}$, we have
\bea
\sum_{\mu\in\pi_{k+1,d}\left(m(\lambda)\right)} c_{\lambda,\mu} e_\mu = \sum_{\mu\in\pi_{k+1},|\mu|\geq d+1, \mu_1\leq M} c_\mu e_\mu.
\nonumber
\eea

For any $(k,r,k+1)$-non-admissible partition $\lambda\in\pi_{k+1}$ such that $\lambda_1\leq M$,
 there exists some $d$ and $m$ so that $\lambda\in\pi_{k+1,d}(m)$.
Moreover, the set $\pi_{k+1,d}(m)$ contains at most one $(k,r,k+1)$-non-admissible partition $\lambda$,
 and for all $\mu \in \pi_{k+1,d}(m)$ such that $\mu \neq \lambda$,
 we have $\mu \succ \lambda$.
Therefore $e_{\lambda}$ can be written in $R_{M,k+1}/{\cal J}_{M,k+1}$ as follows:
\bea
e_{\lambda} = \sum_{\mu \succ \lambda, \mu_1\leq M} c_\mu' e_\mu.
\nonumber
\eea

Let $\lambda \in \pi_{n}$ be a $(k,r,n)$-non-admissible partition such that $\lambda_1\leq M$.
Then there exists $i$ such that $\lambda_i-\lambda_{i+k}<r$.
We set $\mu:=(\lambda_i,\cdots,\lambda_{i+k})\in\pi_{k+1}$
 and $\nu:=(\lambda_1,\cdots,\lambda_{i-1},\lambda_{i+k+1},\cdots,\lambda_n)$.
Since $\mu$ is $(k,r,k+1)$-non-admissible, from the above argument,
 we can rewrite $\mu$ as a linear combination of greater monomials $\{e_{\mu'};\mu'\succ\mu\}$ in $R_{M,k+1}/{\cal J}_{M,k+1}$.
Hence $e_\lambda$ can be rewritten in $R_{M,n}/{\cal J}_{M,n}$ as follows:
\bea
e_\lambda&=&e_\mu e_\nu\nonumber\\
&=&\left(\sum_{\mu'\succ\mu,\mu'_1\leq M} c_{\mu'} e_{\mu'}\right) e_\nu\nonumber \\
&=&\sum_{\lambda'\succ\lambda,\lambda_1'\leq M} c_{\lambda'}e_{\lambda'}.
\nonumber
\eea
Here, in the last $=$, we set $\lambda':=\mu'\cup\nu$.

If $e_{\lambda'}$ is still $(k,r,n)$-non-admissible for some $\lambda'$,
 we further rewrite $e_{\lambda'}$ as a linear combination of greater monomials.
Since $\{ \lambda \in \pi_n; \lambda_1\leq M \}$ is a finite set, this procedure stops in finite times.
\end{proof}

\begin{cor}\label{cor:upper}
$\dim J^{(k,r)}_M\leq \sharp\{\lambda \in \pi_n;\lambda$ is $(k,r,n)$-admissible and $\lambda_1\leq M\}$.
\end{cor}

By Corollary \ref{cor:lower} and Corollary \ref{cor:upper}, we complete the proof of Theorem \ref{thm:main}.

%%%%%%%%%%          section macdonald          %%%%%%%%%%

\section{Application to Macdonald symmetric polynomials}
We can apply the method in Section 3 to a proof of Theorem \ref{thm:old}.

In \cite{Mac}, symmetry relations (Ch. VI, (6.6)) and special values (Ch. VI, (6,11')) of Macdonald symmetric polynomials have been given.
By a combinatorial argument similar to the one employed in this paper,
 we see that for any $(k,r,n)$-admissible partition $\lambda$,
 the multiplicity of the factor $(1-t^{k+1}q^{r-1})$ in r.h.s. of (6,11') is $[\frac{n}{k+1}]$.
Moreover, for $\mu\in\pi_{\eta,N}$ or $\pi_{\eta,N}'$,
 the same results as Lemma \ref{lem:thickwheel} follow as well.
Hence from symmetry relations,
 through the same argument as Theorem \ref{thm:generic},
 we conclude that the Macdonald symmetric polynomial is well-defined and satisfies the wheel conditions if $\lambda$ is $(k,r,n)$-admissible.

~

%%%%%%%%%%          references          %%%%%%%%%%

\end{document}